\ifx\control\undefined

\documentclass[10pt, leqno, twoside]{article}
\usepackage{amssymb, amsmath} 
\usepackage{latexsym} 
\usepackage{graphics} 
\usepackage{url} 
\usepackage{enumerate}
\usepackage[all]{xy}
\usepackage{graphicx}
\usepackage{float}
\usepackage{amsrefs}
\usepackage{color}

\usepackage{authblk}
\usepackage{blindtext}

\usepackage{enumitem}

\numberwithin{equation}{section}

\newtheorem{theo}{Theorem}
\newtheorem{coro}[theo]{Corollary}

\newtheorem{prop}[theo]{Proposition}
\newtheorem{exam}[theo]{Example}

\newtheorem{open}[theo]{Open Question}

\newtheorem{property}[theo]{Property}

\newenvironment{proo}[1][\proofname]{\normalfont{\itshape
#1{:}}\quad\mdseries\ignorespaces}
{{$\Box$}{\vskip\belowdisplayskip}}
\newcommand{\proofname}{Proof}
\newtheorem{defi}[theo]{Definition}
\newtheorem{rem}[theo]{Remark}

\DeclareMathOperator{\modulo}{mod}

\begin{document}

\fi

\title{Generalization of Kimberling's concept of triangle center for other polygons}


\author{L. Felipe Prieto-Mart\'inez\thanks{Department of Mathematics, Universidad Aut\'onoma de Madrid (Spain), luisfelipe.prieto@uam.es} \, and\, Raquel S\'anchez-Cauce\thanks{Department of Artificial Intelligence, Universidad Nacional de Educaci\'on a Distancia (Spain), rsanchez@dia.uned.es}}

\date{\today}

\maketitle

\maketitle

\begin{abstract} In this article we introduce a general definition of the concept of center of an $n$-gon, for $n\geq 3$, generalizing the idea of C. Kimberling for triangle. We define centers associated to functions instead of to geometrical properties. We discuss the definition of those functions depending on both, the vertices of the polygons or the lengths of sides and diagonals. We explore the problem of characterization of regular polygons in terms of these $n$-gon center functions and we study the relation between our general definition of center of a polygon and other approaches arising from Applied Mathematics.

\textbf{Keywords:} Polygon, Triangle Center, Center Function, Center of a Polygon.

\textbf{Mathematics Subject Classification (2010):} Primary 51M04 ,  Secondary 51M15
\end{abstract}

\section{Introduction}

C. Kimberling \cite{K.O, K.FE} in the second half of the 20th century decided to give an unified definition for triangle center, including  the classical centers (incenter, barycenter, circumcenter and orthocenter) and much others (Steiner Point, Fermat point,\ldots). Moreover, he created an encyclopedia \cite{K.E} trying to contain all  known triangle centers. His new idea was to consider triangle centers as functions of the set of lengths, instead of loci.

\medskip

Following this spirit, we provide a definition of center of a polygon as a function ($n$-gon center functions). We also provide a geometric interpretation of this center functions as points in the plane.  Although we still can find some works in the literature studying ``centers'' of polygons (see for example \cite{AS.Q, K.CCP, MM.E}), as far as we know there is no general definition for this concept.

\medskip

The main obstacle is that $n$-gons are not determined by their sidelegths for $n\geq 4$,  so we have defined \emph{centers} as  functions of the vertices (Definition \ref{def.main}). To connect with the definition given by Kimberling for triangles, we also provide an alternative and equivalent definition for the concept of \emph{center} involving, not only sidelengths, but also the lengths of the diagonals, which seems also to be fruitful to encompass some of the well known examples of polygon centers (Definition \ref{def.main2}).

\medskip

In \cite{AS.Q} the authors already studied the problem  of exploring \emph{``the degree of regularity implied by the coincidence of two or more''} centers for quadrilaterals. In that article only  the center of mass of the four vertices,  the center of mass of  the four sides,  the center of mass of the whole figure considered as a lamina of uniform density and the Fermat–Torricelli center are considered. This problem is related to the fact that for squares those four centers coincide. In our new general setting all the centers (in their geometric interpretation as points in the plane) also coincide for  regular polygons (Proposition \ref{prop.reg}). In relation to this, we study the problem of \emph{characterization} of  equiangular, equilateral  and regular polygons by means of one or more  centers (Theorem \ref{th.equiangular}, Theorem \ref{th.equilateral} and Corollary \ref{coro.reg}).

\medskip

We also explore other possibilities for defining \emph{center of a polygon}, that would include some interesting examples (Definition \ref{defi.implicit} and Definition \ref{defi.optimization}) arising from other areas of Mathematics. We also briefly study their relation to our concept of center of a polygon.

\medskip

This article is structured as follows. In  Section \ref{section.tc} we review the main aspects about Kimberling's definition of triangle center. Our definition  of $n$-gon center function and its geometric interpretation appear in Section  \ref{section.defi}, together with a discussion of its suitability: we justify that it satisfies some properties that we would expect for a point to be called \emph{center}.  In Section \ref{section.equiv} we provide an equivalent defintion of center function  in tems of sidelegths and lengths of diagonals, as announce before. Section \ref{section.exam} includes some examples of classic polygon centers that are subsumed in our new definition of \emph{center}. The problem of characterization  of  equiangular, equilateral  and regular polygons by means of one or more  centers is studied in Section \ref{section.characterization}. The other possibilities for a definition of \emph{center} are discussed in Section \ref{section.othercenter}. Finally in Section \ref{section.final} we briefly comment  some topics related to this concept of polygon center and propose some future worklines (in relation to equivariant maps, computational geometry and others). We also include some Open Questions to close up this work.

\section{Basics about triangle centers} \label{section.tc}

As we explained in the introduction, C. Kimberling in his works \cite{K.O, K.FE, K.E} decided to define \emph{triangle centers} as functions (the so-called \emph{triangle center functions}) instead of loci in the plane. These functions also allow a geometric interpretation as  points in the plane, via trilinear coordinates, as it is described below.

\medskip

Unlike what happens for $n$-gons for $n\geq 4$, triangles are determined by their sidelengths:

\begin{rem}\label{rk:identification} Every triangle $T$ can be identified with the tuple of its three sidelengths $(a,b,c)$ (placed in clockwise order), up to congruences. 
\end{rem}

\noindent According to this identification, we denote the set of all triangles as 
$$\mathcal T = \{(a,b,c) \in \mathbb{R}^3_+ \, :\, a +b > c,\ a+c > b, \ b+c >a\} . $$

\medskip

Before proceeding with the definition, let us say that we denote by $[a:b:c]$ the points in the real projective plane $P_2\mathbb R$ with the usual convention $[a:b:c]=[a':b':c']\Leftrightarrow \frac{a}{a'}=\frac{b}{b'}=\frac{c}{c'}$. Then we have:

\begin{defi}[Kimberling's definition of triangle center  \cite{K.O, K.FE, K.E}] \label{defi.TC} A real-valued function $f$ of three real variables $a,b,c$ is a \emph{triangle center function} if it satisfies the following properties:

\begin{enumerate}

\item[(i)] Homogeneity: there exists some constant $n\in\mathbb N$ such that for all $t\in\mathbb R_{\geq 0}$ we have  $f(ta,tb,tc)=t^nf(a,b,c)$.

\item[(ii)] Bisymmetry in the second and third variables: for all $(a,b,c)\in\mathcal T$, we have $f(a,b,c)=f(a,c,b)$.

\end{enumerate}

\noindent Define also the \emph{coordinate map} $\varphi:\mathcal T\to P_2\mathbb R$ given by:
$$\varphi_f(a,b,c)=[f(a,b,c):f(b,c,a):f(c,a,b)] .$$

\end{defi}

$\varphi_f(a,b,c)$ can be also interpreted as trilinear coordinates (the first one corresponding to the side of length $a$, the second one to the side of lenght $b$ and the last one to the one of length $c$). Thanks to these two interpretations, we can think of a center as, actually, a point in the plane.

 \begin{figure}[h]
  \centering
    \includegraphics[width=50mm]{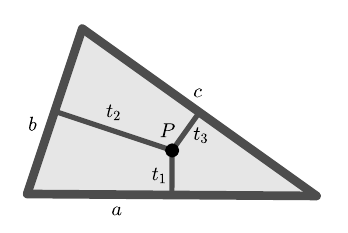}
    \caption{There is a unique point $P$ in the plane such that the relative distances from $P$ to each one of the sides are  $[t_1:t_2:t_3]$. $[t_1:t_2:t_3]$ are then said to be the trilinear coordinates of  $P$.}
\end{figure}

\begin{exam}[circumcenter, see \cite{K.E}] The triangle center function corresponding to the circumcenter is $f(a,b,c)=a(b^2+c^2-a^2)$, since the trilinear coordinates corresponding to the circumcenter are:
$$[t_1:t_2:t_3]=[a(b^2+c^2-a^2):b(c^2+a^2-b^2):c(a^2+b^2-c^2)]. $$

\end{exam}

Kimberling's definition of \emph{triangle center function} ensures that:

\begin{property}[the coordinate map is well defined] \label{property.wd} The trilinear coordinates of the center are associated to the triangle independently of the labelling of the sides (but obviously re-ordered).

\end{property}

\begin{property}[the definition is coherent with respect to similarities] \label{property.cs} Let $f$ be a triangle center function. Let $T_1=(a,b,c)$, $T_2=(a',b',c')$ be two similar triangles such that   $S$ is the similarity $T_1=S(T_2)$. Let $P_1,P_2$ be the points with trilinear coordinates given by  $\varphi_f(a,b,c)$ and $\varphi_f(a',b',c')$, with respect to each of the triangles. Then  $P_1=S(P_2)$.

\end{property}

In this setting, we see that if we have a triangle $(a_1,a_2,a_3)$ and a permutation $\sigma$ of the set $\{1,2,3\}$, then if:
$$\varphi_f(a_1,a_2,a_3)=[t_1:t_2:t_3]=[f(a_1,a_2,a_3):f(a_2,a_3,a_1):
f(a_3,a_1,a_2)], $$

\noindent we have that:
$$\varphi_f(a_{\sigma(1)},a_{\sigma( 2)},a_{\sigma(3)})=[t_{\sigma(1)}:t_{\sigma(2)}:t_{\sigma(3)}] .$$

There is a correspondence between trilinear coordinates and barycentric coordinates. We are interested in this second setting. Kimberling decided to use the first option, although he also claimed the existence of this correspondence:

\begin{rem} If $[t_1:t_2:t_3]$ are the trilinear coordinates of a point $P$ (with respect to the sides of a triangle), then $[at_1:bt_2:ct_3]$ are the barycentric coordinates of this point (with respect to the vertices of this triangle).
\end{rem}

\section{Polygon centers} \label{section.defi}

The concept of triangle center can be generalized to polygons. First we will fix the notation. For us, a polygon is a finite number of straight line segments  connected in a closed chain. We say that the polygon is \emph{non-degenerated} if none of the vertex coincide and that it is \emph{simple} if those segments only intersect with the adjacent element of the chain in a vertex. We will also introduce the following convention, which we will use along the paper:

\begin{rem} \label{rem.conv} Any $n$-gon can be identified with an $n$-uple $(V_1,\ldots,V_n)$ after a labelling of its vertices. This labelling is chosen in such a way that  for $1\leq i\leq n$ the segments  $V_i,V_{(i+1\mod n)}$  are  edges of the polygon, and the rest of segments joining vertices are diagonals.
\end{rem}

Consider the dihedral group $D_n=\langle\rho,\sigma:\rho^n=id,\sigma^2=id,\sigma\rho\sigma=\rho^{-1}\rangle$ (it has $2n$ elements). It can be viewed as a subset of the permutation group of the set $\{1,\ldots,n\}$, determined by
\begin{equation} \label{eq.sigma}\rho(i)=i+1\mod n \, , \qquad\qquad \sigma(i)=2+n-i \mod n. \end{equation}

\noindent But it can also be viewed as a rellabeling (in the sense explained in Remark \ref{rem.conv}) in the set of all $n$-gons:
$$\forall \alpha\in D_n,\qquad \alpha(V_1,\ldots,V_n)=(V_{\alpha(1)},\ldots,V_{\alpha(n)}) .$$

As in the case of triangles, we want to define the \emph{$n$-gon center function} as a function of the vertices $f(V_1,\ldots,V_n)$ (then defined in $(\mathbb R^2)^n$) and the  \emph{coordinate map}. Trilinear coordinates are not a good option to provide a geometric interpretation since they do not extend in a natural way from triangles to $n$-gons for $n\geq 4$. So, we will use barycentric coordinates instead.

\medskip

For a fixed $n$, we denote the set of all $n$-gons by $\mathcal{P}_n\approx (\mathbb R^2)^n$. See that:

\begin{rem}[other domains] \label{rem.domain} Sometimes, we may restrict ourselves to convex $n$-gons (whose vertices satisfy (1) $V_i\neq V_j$ if $i\neq j$ and (2) for any $i$, all vertices (except $V_i,V_{i+1\mod n}$) lie on the same side of the line  defined by $V_i, V_{i+1\mod n}$) or to non-degenerated $n$-gons.
\end{rem}

Now we are ready for generalizing the definition of \emph{triangle center function} by C. Kimberling:

\begin{defi}[main definition of $n$-gon center function] \label{def.main} We say that a real-valued function \break $f(V_1,\ldots,V_n)$ is a \emph{$n$-gon center function} if it satisfies the following properties:

\begin{enumerate}

\item[(1)] Preservation with respect to relabellings: for the symmetry $\sigma \in D_n$:
$$f(V_1,V_2,\ldots,V_n)=f(V_1,V_{\sigma(2)},\ldots,V_{\sigma(n)}) .$$

\item[(2)] Homogeneity: there exists some $k\in\mathbb N$ such that, for all $t\in\mathbb R_{\geq 0}$ we have that $f(tV_1,\ldots,tV_n)=t^kf(V_1,\ldots,V_n)$.

\item[(3)] Preservation with respect to motions: for every rigid motion $T$ in the plane, $f(T(V_1),\ldots,T(V_n))=f(V_1,\ldots,V_n)$.

\end{enumerate}

\noindent Define also the \emph{coordinate map} $\varphi : \mathcal{P}_n \to P_{n-1} \mathbb R$ given by:
$$\varphi_f(V_1,\ldots,V_n)=[f(V_1,\ldots,V_n): f(V_2,\ldots,V_n,V_1):\ldots: f(V_n,V_1,\ldots,V_{n-1})] .$$

\noindent Note that it is not possible to define the coordinate map for $n$-gons such that 
$$f(V_1,\ldots,V_n)=f(V_2,\ldots,V_n,V_1)=\ldots=f(V_n,V_1,\ldots,V_{n-1})=0.$$

\end{defi}

Coordinates $\varphi_f(V_1,\ldots,V_n)$  (when defined) are interpreted as barycentric coordinates with respect to the vertices. Then, the  \emph{geometric interpretation} of the center of a given  $n$-gon $(V_1,\ldots,V_n)$ associated to the $n$-gon center function $f$ is the point:
\begin{equation} \label{eq.gi} P=\widetilde f(V_1,\ldots,V_n)V_1+\ldots+\widetilde f(V_n,\ldots,V_{n-1})V_n \end{equation}

\noindent  where $\widetilde f(V_1,\ldots,V_n),\ldots ,\widetilde f(V_n,\ldots,V_1)$ are the normalized coordinates
\begin{equation} \label{eq.norm}\begin{cases} \varphi_f(V_1,\ldots,V_n)=[\widetilde f(V_1,\ldots,V_n):\ldots: \widetilde f(V_n,V_1,\ldots,V_{n-1})],\\
\widetilde f(V_1,\ldots,V_n)+\ldots+ \widetilde f(V_n,V_1,\ldots,V_{n-1})=1 \end{cases}\end{equation}

We want the coordinate map to satisfy an analogue of Properties \ref{property.wd} and \ref{property.cs}. The next result ensures that, and may help us to clarify the notation and ideas in this paper.

\begin{theo} Definition \ref{def.main} and the geometric interpretation described in \eqref{eq.gi} provide an analogue to  Properties \ref{property.wd} and \ref{property.cs} for $n$-gon center functions, i.e., 

\begin{itemize}

\item[P1] \emph{(the coordinate map is well defined)} The coordinates given by the coordinate map $\varphi_f$ are associated to each  polygon independently of the labelling (but obviously re-ordered).

\item[P2] \emph{(the definition is coherent with respect to similarities)}  Let $f$ be a $n$-gon center function. Let $N_1=(V_1,\ldots,V_n)$ and $N_2=(V_1',\ldots,V_n')$ be two similar  $n$-gons such that $S$ is the similarity $N_1=S(N_2)$. Let $P_1,P_2$ be the points with barycentric coordinates given by $\varphi_f(V_1,\ldots,V_n)$ and $\varphi_f(V_1',\ldots,V_n')$, with respect to each of the $n$-gons. Then $P_1=S(P_2)$.

\end{itemize} 

\end{theo}

\begin{proo}  To prove P1, see that for every $\alpha\in D_n$
$$\varphi_f(\alpha(V_1,\ldots,V_n))=\alpha(\varphi_f(V_1,\ldots,V_n)), $$

\noindent where $\alpha([t_1:\ldots:t_n])=[t_{\alpha(1)}:\ldots:t_{\alpha(n)}]$. 

\medskip

To prove P2  see that any similarity $S$ can be obtained as a composition $T\circ H$ of a rigid motion $T$ and an homotethy $H=\lambda\cdot id$ fixing the origin. Hence, 
\begin{align*}
\varphi_f(S(V_1),\ldots,S(V_n)) & =\lambda^k \varphi_f(T(V_1),\ldots,T(V_n))  =\varphi_f(T(V_1),\ldots,T(V_n))= \\
&=\varphi_f(V_1,\ldots,V_n).
\end{align*}

\hspace{15cm}\end{proo}

We conclude this section with the following result, which states an important property.  Recall that we say that an $n$-gon is \emph{regular} if it is equiangular and equilateral. Regular $n$-gons can be either convex or star.

\begin{prop}\label{prop.reg} 
For any center function $f$, a regular $n$-gon (convex or star) $(V_1,\ldots,V_n)$  satisfies 
$$f(V_1,\ldots,V_n)=f(V_2,\ldots, V_n,V_1)=\ldots=f(V_n,V_1,\ldots,V_{n-1}), $$

\noindent so, if defined, $\varphi_f(V_1,\ldots,V_n)=[1:\ldots:1]$.

\end{prop}

\begin{proo}
Let $(V_1, \ldots, V_n)$ be a regular $n$-gon and $f$ a center function for this $n$-gon. Then, the $n$-gon $(V_i, V_{i+1}, \ldots , V_n, V_1, \ldots , V_{i-1})$ is also  regular and corresponds to a rotation $T_i$ of $\frac{2\pi (i-1)}{n}$rad of $(V_1, \ldots, V_n)$, i.e.,
$$T_i(V_1, \ldots, V_n) = (T_i(V_1), \ldots, T_i(V_n)) = (V_i, V_{i+1}, \ldots , V_n, V_1, \ldots , V_{i-1}). $$

\noindent  Since a rotation is a rigid motion in the plane, by property (3) of Definition \ref{def.main} we have that \break $f(T_i(V_1), \ldots, T_i(V_n))= f(V_1, \ldots, V_n)$. Thus, if defined, 
\begin{align*}
\varphi_f & (V_1,\ldots,V_n)  =[f(V_1,\ldots,V_n): f(V_2,\ldots,V_n,V_1):\ldots: f(V_n,V_1,\ldots,V_{n-1})] \\
& =[f(V_1,\ldots,V_n): f(T_2(V_1), \ldots, T_2(V_n)):\ldots: f(T_n(V_1), \ldots, T_n(V_n))] \\
&= [f(V_1,\ldots,V_n): f(V_1,\ldots,V_n):\ldots: f(V_1,\ldots,V_n)] = [1:1: \ldots :1] .
\end{align*}

\hspace{15cm} \end{proo}

\section{An equivalent definition of polygon center function in terms of lengths} \label{section.equiv}

The definition of center above is, in some sense, not satisfactory. The first reason is that it may be not inmediate to verify condition (3). And the second one is that, traditionally, some of the most useful centers are described in terms of the sidelengths, not in terms of the coordinates of the vertices.

\medskip

We need again to stablish some conventions.
Let $(V_1,\ldots,V_n)$ be an $n$-gon. We will denote by $d_{ij}$  the length of the segment with endpoints $V_i,V_j$. It is obvious that $d_{ij} = d_{ji}$. If $j=i+1\mod n$, then $d_{ij}$ is a sidelength. We will write  $e_{ij}$ instead of $d_{ij}$ when we want to emphasize that we are referring to sidelengths.

\medskip

An $n$-gon is not completely determined by its sidelengths, some of the lenths of the diagonals are required to determine it up to congruence. The set of all the sidelengths and  of the lengths of the  diagonals of an $n$-gon must satisfy some compatibility conditions. For example, consider a quadrilateral with sidelengths $e_{12},e_{23},e_{34},e_{41}$ and diagonals $d_{13},d_{24}$. According to the Cayley-Menger determinant formula for the volume of a 3-dimensional tetrahedron (see \cite{S.NDIM}) we have that:
$$\left|\begin{array}{c c c c c}0 & 1 & 1 & 1 & 1 \\ 1 & 0 & e_{12}^2 & d_{13}^2 & e_{14}^2 \\ 1 & e_{21}^2 & 0 & e_{23}^2 & d_{24}^2 \\ 1 & d_{31}^2 & e_{32}^2 & 0 & e_{34}^2 \\ 1 & e_{41}^2 & d_{42}^2 & e_{43}^2 & 0 \end{array}\right|=0 .$$

So, sometimes we will identify an $n$-gon with the $(n+(n-3))$-uple
$$(e_{12},\ldots,e_{n1},d_{13},\ldots, d_{1,n-1})$$

\noindent  and sometimes, if it is more simple for the corresponding formulas, with the $n\cdot (n-1)$-uple $(d_{ij})$ with $i\neq j$, taking in mind that some of the entries $d_{ij}$ are redundant.

\medskip

In this setting we can define the $n$-gon center function as a real-valued  function  $g$ depending on  the sidelengths and the lengths of the diagonals (defined then in $\mathbb{R}^{(n(n-1))}$), instead of the vertices, as follows.

\begin{defi}[definition of $n$-gon center function in terms of lengths] \label{def.main2} We say that a real-\break valued function $g(d_{ij})$, $i,j= 1, \ldots,n$, is a \emph{$n$-gon center function} if it satisfies the following properties:

\begin{enumerate}
\item[(1')] Preservation with respect to relabellings: for the  symmetry $\sigma \in D_n$:
$$g(d_{ij})=g(d_{\sigma(i),\sigma(j)}). $$

\item[(2')] Homogeneity: there exists some $k\in\mathbb N$ such that, for all $t\in\mathbb R_{\geq 0}$, we have that $g(t\cdot d_{ij})=t^k\cdot g(d_{ij})$.

\end{enumerate}

\noindent Define also the \emph{coordinate map} $\varphi : \mathcal{P}_n \to P_{n-1} \mathbb R$ given by:
$$\varphi_g(d_{ij})=[g(d_{ij}):g(d_{\rho(i),\rho(j)}):\ldots:g(d_{\rho^{n-1}(i),\rho^{n-1}(j)})] ,$$
for $\rho\in D_n$.

\noindent Note that it is not possible to define the coordinate map for $n$-gons such that:
$$g(d_{ij})=g(d_{\rho(i),\rho(j)}))=\ldots=g(d_{\rho^{n-1}(i),\rho^{n-1}(j)})=0. $$

\end{defi}

In this context, the geometric interpretation of the \emph{center} is again (in those cases where the coordinate map is defined):
$$\widetilde g(d_{ij})V_1+\widetilde g(d_{\rho(i),\rho(j)}))V_2+\ldots+\widetilde g(d_{\rho^{n-1}(i),\rho^{n-1}(j)})V_n$$

\noindent where $\widetilde g(d_{ij}),\ldots,\widetilde g(d_{\rho^{n-1}(i),\rho^{n-1}(j)})$ are the normalized coordinates as done in \eqref{eq.norm}.

\medskip

The next result ensures that both definitions of center function are compatible.

\begin{theo}[equivalence between definitions  \ref{def.main} and \ref{def.main2}]  \label{th.equiv} Given a $n$-gon center function \break  $f(V_1,\ldots,V_n)$, it is possible to find an $n$-gon center function $g(d_{ij})$  such that the geometric interpretations (if they exist) of the centers corresponding to $f$ and $g$ coincide for every element in $dom(f)\subset \mathcal P_n$, and viceversa.

\end{theo}

\begin{proo}
Supose that we have an $n$-gon with vertices  $(V_1,\ldots,V_n)$. The vertices determine univocally the lengths $d_{ij}=\left|V_i-V_j\right|$.  On the other hand, the lengths $d_{ij}$ determine modulo congruence the vertices $(V_1,\ldots,V_n)$. So we can consider the vertices $V_k$ as  functions of the lengths $V_k(d_{ij})$, if we impose $V_1=(0,0)$, $V_2=(d_{12},0)$ and that the rest of the vertices are ordered clockwise.
\medskip

First, see that if $f(V_1,\ldots,V_n)$ is a center function in the sense of Definition \ref{def.main}, then it is easy to find an associated center function $g(d_{ij})$ in the sense of Definition \ref{def.main2} of the form:
$$g(d_{ij})=f(V_1(d_{ij}),\ldots,V_n(d_{ij})) .$$

We just have to prove that if $f$ satisfies properties (1), (2) and (3), then $g$ satisfies properties (1') and (2'). For $\sigma\in D_n$ as defined in  \eqref{eq.sigma} we have:
\begin{align*}
 g(d_{\sigma(i),\sigma(j)})&=f(V_1(d_{\sigma(i),\sigma(j)}),V_{\sigma(2)}(d_{\sigma(i),\sigma(j)}),\ldots,V_{\sigma(n)}(d_{\sigma(i),\sigma(j)})) =\\
&=f(V_1(d_{ij}),\ldots,V_n(d_{ij})) = g(d_{ij}), \end{align*}

\noindent by properties (1) and (3). So, property (1') holds. On the other hand:
\begin{align*}
g(t\cdot d_{ij})& =f(V_1(t\cdot d_{ij}),\ldots,V_n(t\cdot d_{ij}))=f(t\cdot V_1(d_{ij}),\ldots,t\cdot V_n(d_{ij}))= \\
& =t^n\cdot f(V_1(d_{ij}),\ldots,t_n(V_n(d_{ij})))=t^n\cdot g(d_{ij}) ,
\end{align*}
and so property (2') also holds.

\medskip

Next, see that if $g(d_{ij})$ is a center function in the sense of Definition \ref{def.main2}, then it is easy to find an associated center function $f(V_1,\ldots,V_n)$ in the sense of Definition \ref{def.main} of the form:
$$f(V_1,\ldots,V_n)=g(\left|V_i-V_j\right|) .$$

\noindent We just have to prove that if $g$ satisfies properties (1') and (2'), then this $f$ satisfies properties (1), (2), and (3). First see that, for $\sigma \in D_n$ as defined in \eqref{eq.sigma} we have that:
$$f(V_1,V_{\sigma(2)},\ldots,V_{\sigma(n)})=g(\left|V_{\sigma(i)}-V_{\sigma(j)}\right|)=g(\left|V_i-V_j\right|)=f(V_1,\ldots,V_n) ,$$
by property (1'). Hence, property (1) holds. Now, see that:
\begin{align*} f(t\cdot V_1,\ldots,t\cdot V_n)&=g(\left|t\cdot V_i-t\cdot V_j\right|)=g(t\cdot \left|V_i-V_j\right|)=t^n\cdot g(\left|V_i-V_j\right|)\\
&=t^n\cdot f(V_1,\ldots,V_n) ,\end{align*}
by property (2'). So, property (2) also holds. Finally, the proof of property (3) is inmediate: all congruent $n$-gons have the same sidelengths and diagonals. \end{proo}


\section{Some examples} \label{section.exam}

In this section we present some of the more relevant centers for polygons. Most of them arise from important problems in Applied Mathematics, and can be naturally defined as an affine combination of the vertices. The coefficients of this combinatio are functions of either the vertices or the sidelengths and lengths of the diagonals. Some of those examples can be found in \cite{AS.Q} in the particular case $n=4$ (quadrilaterals).

\begin{exam}[centroid, barycenter or center of mass of the vertices] The barycenter of a polygon with vertices $V_1,\ldots,V_n$, or the center of mass of the vertices (provided that all the vertices  have the same weight) is the point:
\begin{equation} \label{eq.bary} G_0(V_1,\ldots,V_n)=\frac{1}{n}V_1+\ldots+\frac{1}{n}V_n . \end{equation}

\noindent So, the associated center function can be chosen to be $f_0(V_1,\ldots,V_n)=1$  and the coordinate map is $\varphi_{f_0}(V_1, \ldots, V_n) = [1:\ldots :1]$ (recall that the coefficients for the affine combination are the normalized ones).
\end{exam}

\begin{exam}[center of mass of the perimeter of a convex polygon]
The center of mass of the perimeter of a convex polygon (provided that all the points in the perimeter have the same weight) with vertices $V_1,\ldots,V_n$  is the point (see \cite{K.CCP}):
$$ G_1(V_1,\ldots,V_n) = \sum_{i=1}^n \dfrac{e_{i-1(\modulo n),i}+e_{i,i+1(\modulo n)} }{2 \sum_{j=1}^n e_{j,j+1 (\modulo n)} } \cdot V_i . $$
So, the associated center function can be chosen as 
$$g_1(d_{ij})= {e_{n1}+e_{12} }$$

\noindent  and the coordinate map is  $\varphi_{g_1}(V_1, \ldots, V_n) =[(e_{n1}+e_{12}):\ldots:(e_{n-1,n}+e_{n1})]$.
\end{exam}

\begin{exam}[centroid of the polygonal lamina] The centroid of a polygonal lamina with vertices $V_1,\ldots,V_n$ is the point (see \cite{B.CPL}):
\begin{equation} \label{eq.exmalo} G_2(V_1,\ldots,V_n) = \sum_{i=1}^n\left(\frac{ V_{i-1 (\modulo n)}\wedge V_i+V_i\wedge V_{i+1 (\modulo n)}}{3\sum_{j=1}^n V_j\wedge V_{j+1 (\modulo n)}}\right) \cdot V_i ,\end{equation}
\noindent where $V_m\wedge V_{k} = x_m y_k - x_k y_m$, for $V_m = (x_m, y_m)$ and $V_k = (x_k, y_k)$. This is not an affine combination since
$$\sum_{i=1}^n\left(\frac{ V_{i-1 (\modulo n)}\wedge V_i+V_i\wedge V_{i+1 (\modulo n)}}{3\sum_{j=1}^n V_j\wedge V_{j+1 (\modulo n)}}\right)=\frac{2}{3}\neq 1,$$

\noindent so this is not the geometric interpretation of a center in our sense.

\medskip

To include this important center in the setting of our definition we are going to modify \eqref{eq.exmalo}. $G_2(V_1,\ldots,V_n)$, if the $n$-gon $(V_1,\ldots,V_n)$ is convex, can be computed via ``geometric decomposition'' as
$$\frac{\sum_{i=1}^n (G_2(V_i,V_{i+1 (\modulo n)},B))\cdot AREA(V_i,V_{i+1 (\modulo n)},B)}{\text{AREA}(V_1,\ldots,V_n)}, $$

\noindent for $B=G_0(V_1,\ldots,V_n)$ (see \eqref{eq.bary}),  which, according to the \emph{Shoelace Formula} (see \cite{P.SF}) and to the formula of the centroif of a triangular lamina (the classical centroid of the triangle), equals to
$$\frac{\sum_{i=1}^n (\frac{1}{3}V_i+\frac{1}{3}V_{i+1 (\modulo n)}+\frac{1}{3}B)\cdot \frac{1}{2}\left|(B-V_i)\wedge (B-V_{i+1 (\modulo n)})\right|}{\frac{1}{2}\sum_{i=1}^n\left|(B-V_i)\wedge (B-V_{i+1 (\modulo n)})\right|} =$$
\begin{equation}\label{eq.exbueno}=\frac{\sum_{i=1}^n(C_1(i)+C_2(i)+C_3(i))}{3\sum_{i=1}^n \left|(B-V_i)\wedge(B-V_{i+1})\right|}\end{equation}

\noindent where:
\begin{align*} C_1(i)&= \left|(B-V_i)\wedge(B-V_{i+1(\modulo n)})\right|,\\
 C_2(i)&=\left|(B-V_{i-1 (\modulo n)})\wedge (B-V_i)\right|,\\
C_3(i)&=\frac{1}{n}\sum_{j=1}^n\left|(B-V_j)\wedge(B-V_{j+1 (\modulo n)})\right|.\end{align*}

\noindent Expression \eqref{eq.exbueno} does correspond to the geometric interpretation of a center with center function
\begin{align*} f_2(V_1,\ldots,V_n)=&\left|(B-V_1)\wedge(B-V_{2})\right|+\left|(B-V_{n})\wedge (B-V_1)\right|+\\
&+\frac{1}{n}\sum_{j=1}^n\left|(B-V_j)\wedge(B-V_{j+1 (\modulo n)})\right|.\end{align*}

\end{exam}

\begin{exam}[medoid] \label{ex:medoid}
The medoid of the set of vertices $V_1,\ldots,V_n$ is the point $G_3$ such that (see, for example, the recent work \cite{B.M}):
$$G_3 ={\arg \min}_{V \in\{V_1,\ldots,V_n\}}\sum_{i=1}^n \left|V-V_i \right| .$$

\noindent The medoid is not well defined for any $n$-gon (this minimum may be reached by two or more of the vertices). The medoid can also be considered as a center in our sense. In this case the center function is:
$$f_3(V_1,\ldots,V_n)=\begin{cases}1 & \text{if }V_1=\min_{V\in\{V_1,\ldots,V_n\}}\sum_{i=1}^n \left|V-V_i \right| , \\ 0 & \text{in other case.} \end{cases}$$

\end{exam}

\section{Characterization of $n$-gons using centers (specially mention to quadrilaterals)} \label{section.characterization}

The idea of characterizing regular polygons  using $n$-gon  center functions was one of the main reasons of our interest in this topic, in connection to other geometric problems. This study was already started in \cite{AS.Q} for quadrilaterals.  We say that:

\begin{defi}
A set of center functions $f_1,\ldots,f_k$ with associated coordinate maps  $\varphi_{f_1},\ldots,\varphi_{f_k}$ characterizes a family $\mathcal F$ of $n$-gons if  $(V_1,\ldots,V_n)\in\mathcal F$ if and only if:
$$\begin{cases} f_1(V_1,\ldots,V_n)=f_1(V_2,\ldots,V_n,V_{1})=\ldots=f_1(V_{n},V_1,\ldots,V_{n-1}) ,\\
\vdots \\
f_k(V_1,\ldots,V_n)=f_k(V_2,\ldots,V_n,V_{1})=\ldots=f_k(V_{n},V_1,\ldots,V_{n-1}) . \end{cases}$$
\end{defi}

If a family $\mathcal F$ is characterized by a set of center functions then it must be closed under congruences and  it must contain regular $n$-gons (convex and star, see Proposition \ref{prop.reg}).

\medskip

Regular triangles (for triangles equilaterality and equiangularity are equivalent properties) are characterized  by just one center function. Take for example: $f(V_1,V_2,V_3)=\|V_2-V_3\|$.  Equiangular quadrilaterals (rectangles and their non-simple version called \emph{crossed rectangles}), provided that they are non-degenerated, are also characterized by one center function:
$$f (V_1,V_2,V_3, V_4) = \frac{\vec{V_1V_2}\cdot \vec{V_1V_4}}{\|\vec{V_1V_2}\|\cdot \|\vec{V_1V_4}\|} .$$

This is not so trivial: the cosine of two angles being equal does not imply the angles are equal but complementary. But in this case this is not a problem since the sum of the angles of a non-degenerated quadrilateral must be less or equal to $2\pi$rad.

However, there is no center or set of centers characterizing either equilateral quadrilateral (rhombi) or regular quadrilaterals (squares). The following results formalize this idea:

\begin{theo} \label{th.equiangular} Equiangular $n$-gons can be characterized by one center function, provided that they are non-degenerated and convex.

\end{theo}

\begin{proo} The $n$-gon center function that characterizes equiangular $n$-gons (provided that they are convex and so the angle between two adjacent sides is less than $\pi$ rad) is again
\begin{equation} \label{eq.centercos} f_1 (V_1,\ldots, V_n) = \frac{\vec{V_1V_2}\cdot \vec{V_1V_n}}{\|\vec{V_1V_2}\|\cdot \|\vec{V_1V_n}\|}.\end{equation}

\hspace{15cm}\end{proo}

\begin{figure}[h]
  \centering
    \includegraphics[width=50mm]{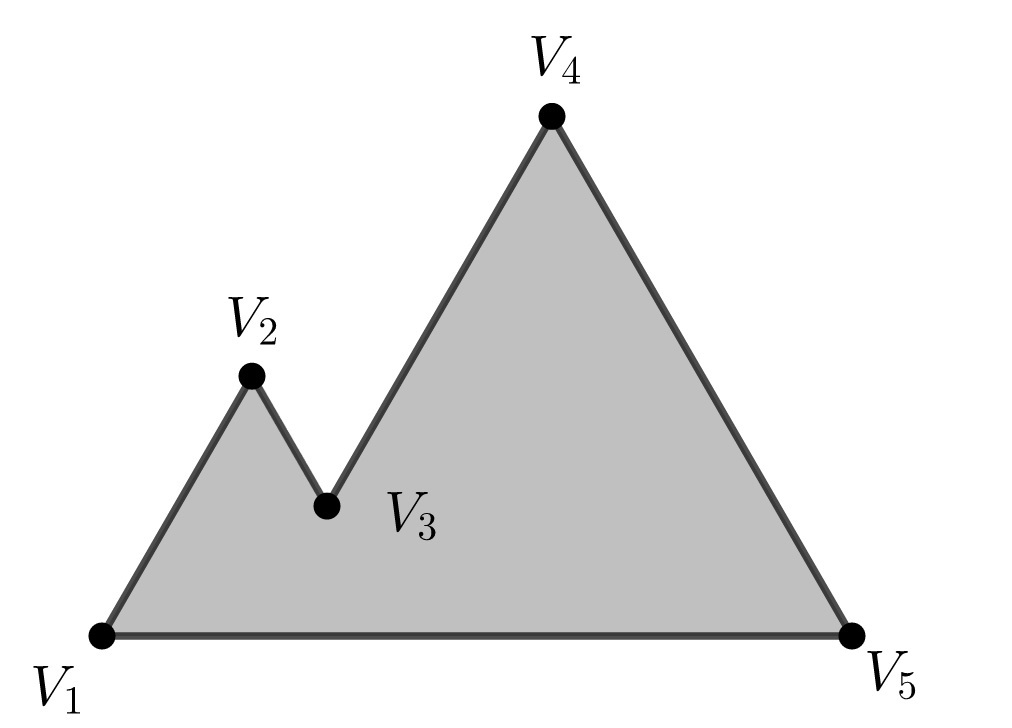}
   \caption{This pentagon is obtained as the union of two equilateral triangles. So all the angles equal to $\pi/3$, except the one over $V_3$, which equals $2\pi-\pi/6$.}
\label{fig.mountain}
\end{figure}

In Figure \ref{fig.mountain} we show a pentagon which, despite not being equiangular, is also included in a family of polygons characterized by the center \eqref{eq.centercos}.

\begin{theo} \label{th.equilateral}

For $n\geq 3$ being an odd number, equilateral $n$-gons can be characterized by one center function.

For $n\geq 4$ being an even number, equilateral $n$-gons cannot be characterized using $n$-gon center functions.  This is a consequence of the fact that any $n$-gon center function $f$ must satisfy:
$$(V_1,\ldots,V_n)\in \mathcal S_n \Rightarrow  f(V_1,\ldots,V_n)=\ldots=f(V_n,V_1,\ldots,V_{n-1}) ,$$

\noindent for $\mathcal S_n$ being the family of equiangular $n$-gons $(V_1,\ldots,V_n)$ such that $\|V_\rho(i)-V_\rho(j)\|=\|V_{\sigma(i)}-V_{\sigma(j)}\|$ (generalization of rectangles for $n\geq 6$).  This family can be characterized by one $n$-gon center function.

\end{theo}

\begin{proo}  The $n$-gon center function that characterizes  equilateral $n$-gon for $n$ odd is
$$f_2(V_1,\ldots,V_n)=\|V_{(n+1)/2}-V_{(n+3)/2}\| .$$

The fact that any $n$-gon center function must characterize $\mathcal S_n$ is inmediate from the ``preservation with respect to relabellings property''. The center function that characterizes $\mathcal S_n$ is:
$$f_3(V_1,\ldots,V_n)=\|V_{n/2}-V_{(n/2)+2}\| . $$

\hspace{15cm} \end{proo}

Theorems \ref{th.equiangular} and \ref{th.equilateral} together imply:

\begin{coro}\label{coro.reg}  For $n\geq 3$ odd, regular $n$-gons  can be characterized by two center functions, provided that they are convex. For $n\geq 4$ even, they cannot.

\end{coro}

See that Theorems 3.1, 3.2, 3.3, 3.4, 3.5 in \cite{AS.Q} are compatible with the results proved here, although the authors are only interested in some particular quadrilateral centers.

\section{Other possible definitions of center} \label{section.othercenter}

Some of the centers arising from Applied Mathematics are defined by an implicit equation involving the vertices, or as solution of a problem of optimization (see Example \ref{ex:medoid}). Besides, the term ``center'' appears in a different setting for compact length spaces in  \cite{MM.E}, from a totally different approach. This leads to the two following alternative definitions of center: 

\begin{defi}\label{defi.implicit} Let $F(P,V_1,\ldots,V_n)$ be a map satisfying the following properties:

\begin{itemize}

\item[(a)] For every $(V_1,\ldots,V_n)$ in $\mathcal P_n$ (additionally non-degeneration property can be required) \break $F(P,V_1,\ldots,V_n)=0$ defines univocally $P$.

\item[(b)] Preservation with repect to labellings: if $F(P,V_1,\ldots,V_n)=0$, then for any element $\alpha\in D_n$ we have
$$F(P,V_{\alpha(1)},\ldots,V_{\alpha(n)})=0 .$$

\item[(c)] Homogeneity: for every $\lambda\in\mathbb R_{\geq 0}$,
$$F(P,V_1,\ldots,V_n)=0\Rightarrow F(\lambda P,\lambda V_1,\ldots,\lambda V_n)=0.$$

\item[(d)] Preservation by motions: for every rigid motion $T$ in the plane,
$$F(P,V_1,\ldots,V_n)=0\Rightarrow F(T(P),T(V_1),\ldots,T(V_n))=0.$$

\end{itemize}

\noindent We say that the point $P$ ensured by (a) is an \emph{implicit center} of $(V_1,\ldots,V_n)$.

\end{defi}

\begin{defi} \label{defi.optimization} Let $G(P,V_1,\ldots,V_n)$ be a real function defined in $\mathcal P_n$  (additionally non-degeration property can be required for the domain) such that:

\begin{itemize}

\item[(a')] For every $(V_1,\ldots,V_n)$ in $\mathcal P_n$ there exists a unique
$$X(V_1,\ldots,V_n)={\arg\min}_{P\in\mathbb R^2}(G(P,V_1,\ldots,V_n)).$$

\item[(b')] Preservation with repect to labellings: for any element $\alpha\in D_n$ we have
$${\arg\min}_{P\in\mathbb R^2}(G(P,V_1,\ldots,V_n))={\arg\min}_{P\in\mathbb R^2}(G(P,V_{\alpha(1)},\ldots,V_{\alpha(n)})).$$

\item[(c')] Homogeneity: for every $\lambda\in\mathbb R_{\geq 0}$,
$$\lambda\cdot{\arg\min}_{P\in\mathbb R^2}(G(P,V_1,\ldots,V_n))= {\arg\min}_{P\in\mathbb R^2}(G(P,\lambda V_{1},\ldots,\lambda V_{n})).$$

\item[(d')] Preservation by motions: for every rigid motion $T$ in the plane,
$$T\left({\arg\min}_{P\in\mathbb R^2}(G(P,V_1,\ldots,V_n))\right)= {\arg\min}_{P\in\mathbb R^2}(G(P,T(V_{1}),\ldots,T(V_{n}))).$$

\end{itemize}

\noindent We will say that the point $X(V_1,\ldots,V_n)$ ensured by (a') is a \emph{minimal center} of $(V_1,\ldots,V_n)$.
\end{defi}

We have that:

\begin{theo} 
Any implicit center in the sense of Definition \ref{defi.implicit} is a minimal center in the sense of Definition \ref{defi.optimization} and viceversa. Moreover, any center in the sense of Definition \ref{def.main} is an implicit center in the sense of Definition \ref{defi.implicit} (and therefore, a minimal center in the sense of Definition \ref{defi.optimization}).
\end{theo}

\begin{proo} To show that definitions \ref{defi.implicit} and \ref{defi.optimization} are equivalent take 
$$F(P,V_1,\ldots,V_n)=G(P,V_1,\ldots,V_n)-X(V_1,\ldots,V_n) $$

\noindent and, for $X$ the unique point satisfying $F(X,V_1,\ldots,V_n)=0$,
$$G(P,V_1,\ldots,V_n)=\min_{P\in\mathbb R^2}(\text{dist}(P,X)) .$$

 The proof of the second statement is inmediate taking 
$$F(P,V_1,\ldots,V_n)=(\widetilde f(V_1,\ldots,V_n)V_1+\ldots+\widetilde f(V_n,V_1,\ldots,V_{n-1})V_n)-P . $$

\hspace{15cm} \end{proo}

Some examples of well-known points usually called ``centers'' that could be naturally included in this different definitions could be:

\begin{exam}[geometric median of the vertices] \label{ex.gm}
The geometric median of the set of vertices \break $V_1,\ldots,V_n$ of an $n$-gon is the point $X$ minimizing the sum of distances to the vertices. Thus, it could be naturally considered as a \emph{minimal center} defined by (see \cite{D.GM}):
$$ X= {\arg \min}_{P \in \mathbb{R}^2} \sum_{j=1}^n \left|V_j-P\right| . $$

\noindent  Provided that $X$ is distinct from any vertex, it can be also described as  an \emph{implicit center} by the formula:
$$\sum_{j=1}^n\frac{V_j-X}{\left|V_j-X\right|} =0 .$$

\noindent It is known that there is no explicit ``simple'' formula for $G$ or its coordinates (see \cite{B.NP}).

\end{exam}

\begin{exam}[Chebyshev center]  The Chebyshev center of a bounded set $Q$ is the center $X$ of the minimal-radius ball enclosing the entire set $Q$ (see \cite{B.CC}). It is described as a minimal center by the formula:
$$ X= {\arg \min}_{P \in \mathbb{R}^2} \left (\max_{V \in Q} \|V - P\|^2 \right ) . $$

\end{exam}

\section{Final comments} \label{section.final}

During the development of this article, some questions have arisen:

\begin{open} Can regular $n$-gons, for $n$-odd, be characterized by only one center function?

\end{open}

\begin{open} What de we know about the Characterization Problem when we do not have the restriction of the polygons being convex?

\end{open}

\begin{open} Is any \emph{implicit  center} in the sense of Definition \ref{defi.implicit} a \emph{center} in the sense of Definition \ref{def.main}? See that  Example \ref{ex.gm} shows that the corresponding \emph{center function} may not be trivial at all to find. 
\end{open}

$\mathcal P_n$ is naturally a  $D_n$-space (a topological space endowed with a group of symmetries, see \cite{M.BU}). In this context, coordinate maps can be understood as $G$-maps. It could be interesting to explore this point of view. In particular, this may connect $n$-gon centers with interesting problems in Plane Geometry such as as the Square Peg Problem and its variants \cite{M.S}.

\medskip

Finally we would like to remark that the study of centers for $k$-dimensional polyhedra ($k\geq 3$, but specially $k=3$) would be of great interest in different areas  (computational geometry and computer vision, for instance), and is a problem still to be explored.

\section*{Acknowledgements}

The second author is supported by a postdoctoral grant (PEJD-2018-POST/TIC-9490) from UNED, co-financed by the Regional Government of Madrid with funds from the Youth Employment Initiative (YEI) of the European Union.

%
%



\end{document}